\documentclass[review]{elsarticle}

\usepackage{epstopdf}
\usepackage{amsmath,amssymb,amsfonts}
\usepackage{algorithmic}
\usepackage{graphicx}
\usepackage{textcomp}
\usepackage{subfigure}
\usepackage{float}
\usepackage{color}
\usepackage{algorithm}
\usepackage{algorithmic}
\usepackage{enumerate}
\usepackage{booktabs}
\usepackage{multirow}
\usepackage{moreverb}
\graphicspath{{./graphics/}}
\newtheorem{theorem}{Theorem}[section]
\newtheorem{lemma}[theorem]{Lemma}
\newtheorem{remark}[theorem]{Remark}

\begin{document}

\begin{frontmatter}

\title{Nonlinear model reference adaptive control approach for governance of the commons in a feedback-evolving game }

\author[mymainaddress]{Fang Yan}

%
%
%

\author[mymainaddress]{Xiaorong Hou\corref{mycorrespondingauthor}}
\cortext[mycorrespondingauthor]{Corresponding author}
\ead{houxr@uestc.edu.cn}

\author[mymainaddress]{Tingting Tian}

\author[address]{Xiaojie Chen\corref{mycorrespondingauthor}}
\cortext[mycorrespondingauthor]{Corresponding author}
\ead{xiaojiechen@uestc.edu.cn}

\address[mymainaddress]{School of Automation Engineering, University of Electronic Science and Technology of China, Chengdu, Sichuan 611731, China}
\address[address]{School of Mathematical Sciences, University of Electronic Science and Technology of China, Chengdu, Sichuan 611731, China}

\begin{abstract}
The governance of common-pool resources has vital importance for sustainability. However, in the realistic management systems of common-pool resources, the institutions do not necessarily execute the management policies completely, which will induce that the real implementation intensity is uncertain. In this paper, we consider a feedback-evolving game model with the inspection for investigating the management of renewable resource and assume that there exists the implementation uncertainty of inspection. Furthermore, we use the nonlinear model reference adaptive control approach to handle this uncertainty. We accordingly design a protocol, which is an update law of adjusting the institutional inspection intensity. We obtain a sufficient condition under which the update law can drive the actual system to reach the expected outcome. In addition, we provide several numerical examples, which can confirm our theoretical results. Our work presents a novel approach to address the implementation uncertainty in the feedback-evolving games and thus our results can be helpful for effectively managing the common-pool resources in the human social systems.
\end{abstract}


\end{frontmatter}


\section{Introduction}

The governance of common-pool resources has vital importance for sustainability \cite{ostrom_90,kraak_ff11,adams_03}.
Cooperation among individuals plays an important role for the governance of common-pool resources. However, how to promote cooperation in this scenario is a huge challenge since rational individuals do not tend to choose to cooperate with others \cite{hardin_68}.
Evolutionary game theory is a powerful theoretical approach for studying under what conditions cooperation can emerge \cite{hofbauer_98,smith_82,weibull_97,nowak_04,sigmund_n10,vasconcelos_ncc13,Chen-JRSI-15,Chen-PRSB-19,Duong-PRSA-21,Jusup-PR-22,Quan-Chaos-23,Chen-ND-23}. In particular,
feedback-evolving game models provides an appropriate framework to study the subtle interdependence of resource and social cooperation for the governance of common-pool resources.

Recently, the study of feedback-evolving game models has received considerable attention \cite{tavoni_jtb12,Chen-SR-14,weitz_pnas16,sugiarto_2017,chen_xj_pcb18,szolnoki_18,shao_yx_epl19,hauert_jtb19,lin_yh_prl19,Quan-JSM-19,wang_x_rspa20,tilman_nc20,Hu-CSF-20,yan_21,cao_lx_21,Shu-PRSA-22,Ma-AMC-23,Liu-elife-23}. On the one hand, the co-evolutionary dynamics of individual actions and resources are studied via feedback-evolving games.
For example, Weitz \emph{et al.} investigated a system in which the payoffs favor unilateral
defection and cooperation, given replete and depleted environments,
respectively and observed oscillations of strategies and the environment in the feedback-evolving game model \cite{weitz_pnas16}. In addition, similar periodic state was reported in asymmetric games due to environmental heterogeneity \cite{hauert_jtb19}.
On the other hand, because the co-evolutionary outcome of strategies and the environment is not always what humans expected, some work further demonstrated that institutions are needed for making the co-evolutionary system reach the desired outcome \cite{tavoni_jtb12,chen_xj_pcb18,yan_21}. For instance, Tavoni \emph{et al.} introduced ostracism in a feedback-evolving game and proved that introducing ostracism can maintain cooperation in resource usage under variable social and environmental conditions \cite{tavoni_jtb12}. Chen and Szolnoki considered the natural growth rate of renewable resources in a feedback-evolving game and found that the desired outcome depends not only on the natural growth rate of resources, but also on the intensity of the institutional inspection and punishment \cite{chen_xj_pcb18}.
Yan \emph{et al.} further considered the delay factor of cooperators' contributions to environmental change in a feedback-evolving game and found that the time delay can induce periodic oscillatory dynamics of cooperation level and environment \cite{yan_21}.

However, it should be noted that there often exists uncertainty of implementation in the realistic management systems of common-pool resource because the institutions do not necessarily execute completely or take the corresponding actions according to the management policies or associated regulations. Hence the real implementation intensity is uncertain or unclear. Such kind of uncertainty makes the existing management regulations often fail to achieve the desired effect on the management of resources \cite{rice_96,dichmont_06,fulton_11,arlinghaus_17}.
Therefore, it is more appropriate to consider a feedback-evolving game with the implementation uncertainty for such co-evolutionary systems.
Furthermore, it is a challenge to make the feedback-evolving game system with the uncertainty factor to reach the desired outcome. Indeed model reference adaptive control (MRAC), as an effective method to handle systems with uncertainty, has attracted much attention in recent years \cite{lavretsky09,yucelen13,kersting18,oliveira18,nguyen18,song19,qu19,yan_22,yan_23}.
In particular, it has been applied to social governance \cite{yermekbayeva_18,yuan_19} and has been proven to be effective for handling the uncertainty of the parameters in predator-prey system. However, it is not clear whether it can handle the feedback-evolving game system with implementation uncertainty.

In this paper, we thereby consider a feedback-evolving game with implementation uncertainty for the co-evolutionary system of resource and strategic behaviors. To be specific, we adopt the feedback-evolving game model in Ref.~\cite{chen_xj_pcb18} where the institutional inspection intensity and punishment are considered for the governance of renewable common-pool resource. Furthermore, we consider that the institutional inspection intensity is uncertain since the inspection activity may not be implemented as required in realistic situations. We then employ the nonlinear MRAC approach to handle this uncertainty in the feedback-evolving game we considered \cite{nam88,scarritt08,asiain21}.
The reference model here we use can produce the same nonlinear feedback-evolving game dynamics as the actual system and its inspection intensity can reach to the preset value. We accordingly present a protocol, which is an update law of adjusting the institutional inspection intensity during a time unit. By means of Lyapunov stability theory, we obtain a sufficient condition under which the update law can drive the actual system to reach the desired outcome. Finally, we provide several numerical examples to confirm our theoretical results.


\section{Model and Method}

We consider the feedback-evolving game model in Ref.~{\cite{chen_xj_pcb18}} where $N$ individuals in a population compete for limited, but partly renewable common-pool resource using two basic consuming strategies, i.e., cooperation and defection. To utilize the resource properly, individuals are allocated a legal amount resource which is assumed as
$b_{\rm l}=\frac{b_{\rm m}y}{R_{\rm m}}$, where $b_{\rm m}$ is the maximal resource portion that each individual is allowed to use per unit of time when the amount of the common pool resource $y(t)$ reaches the carrying capacity of resource pool $R_{\rm m}$.
Cooperators follow the allocation rule and restrain their use to $b_{\rm l}$ amount from the resource.
Defectors ignores the allocation rule and utilizes the common pool more intensively by
getting a larger $b_{\rm v}>b_{\rm l}$ amount satisfying  $b_{\rm v}=b_{\rm l}(1+\alpha)$, where $\alpha>0$  characterizes
the severity of defection. To avoid resource exploitation, consumers are monitored and defection is punished by a centrally organized management.
And the institutional monitored probability during a time unit and the fine on the involved defector are set as $p^{*}$ ($0<p^{*}<1$) and $\beta$ ($\beta>0$), respectively.
Since the inspection activity may not be implemented as required, we suppose that the
institutional monitored probability during a time unit denotes as $\hat{p}$ ($0<\hat{p}<1$), which is uncertain. Given the uncertainty in institutional monitored probability during a time unit, the feedback-evolving game model with this uncertain parameter $\hat{p}$ can be expressed by
\begin{equation}
\left\{
\begin{array}{ll}
\dot{x}(t)=x(t)[1-x(t)][\hat{p}\beta-\frac{\alpha b_{\rm m}y(t)}{R_{\rm m}}]\\
\dot{y}(t)=ry(t)[1-\frac{y(t)}{R_{\rm m}}]-N\frac{y(t)}{R_{\rm m}}b_{\rm m}[1+(1-x)\alpha]\,\,,\label{1.2}
\end{array}
\right.
\end{equation}
where $x(t)$ is the fraction of cooperators in the population at time $t$ and $r$ is the intrinsic growth rate.

\begin{figure}
\centering
\includegraphics[width=5in]{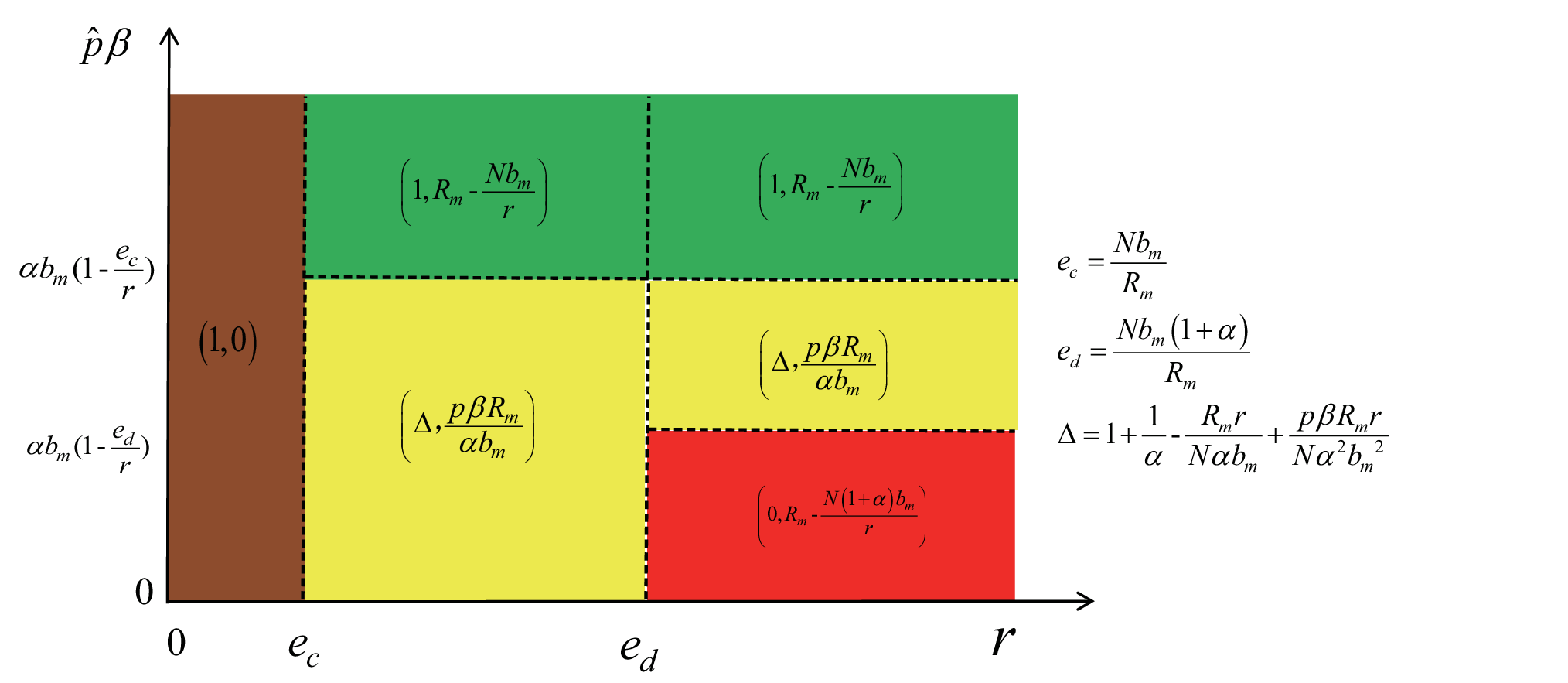}
\caption{A representative plot of evolutionary outcomes on the phase plane. Different colors are used to distinguish qualitatively different solutions in
the parameter space $(r,\hat{p}\beta)$. This plot highlights that the inner dynamical feature of renewable resource could be a decisive factor that can derogate the
expected consequence of punishment. Source: reprinted figure from Ref.~\cite{chen_xj_pcb18}.
}
\label{fig1}
\end{figure}

Through analysis we can know that this equation system has at most five fixed points
which are $(0,0)$, $(1,0)$, $(0,R_{\rm m}-\frac{Nb_{\rm m}(1+\alpha)}{r})$,
$(1,R_{\rm m}-\frac{Nb_{\rm m}}{r})$, and
$(\Delta,\frac{R_{\rm m}\hat{p}\beta}{\alpha b_{\rm m}})$, respectively,
where $\Delta=1+\frac{1}{\alpha}-\frac{rR_{\rm m}}{\alpha b_{\rm m}N}+\frac{\hat{p}\beta R_{\rm m}r}{N\alpha^{2}b_{\rm m}^{2}}$.
For further analysis let us denote $e_{c}=\frac{Nb_{\rm m}}{R_{\rm m}}$ and $e_{d}=\frac{Nb_{\rm m}(1+\alpha)}{R_{\rm m}}$, respectively, representing the gain rates of cooperators and defectors in a population from the common resource. A representative plot of evolutionary outcomes on the phase plane is shown in Fig.~1. We can see that $0<r<e_{c}$ means the resource pool grows slowly, $e_{c}<r<e_{d}$ means the resource pool grows moderately, and $e_{d}<r$ means the resource pool grows fast. Here we are interested in the equilibrium point, $(1,R_{\rm m}-\frac{Nb_{\rm m}}{r})$, which denotes the desired outcome of the co-evolutionary system.
This desired outcome means that all participants choose the cooperation strategy and the renewable resource can be maintained at a sustainable level.
The desired outcome can be obtained when the intrinsic growth rate and the institutional inspection probability during a time unit satisfy these conditions $r>e_{c}$ and $\hat{p}>(\frac{1}{\beta}-\frac{e_{c}}{r\beta})\alpha b_{\rm m}$, respectively.

Note that as long as the institutional inspection probability during a time unit satisfies the condition $\hat{p}>(\frac{1}{\beta}-\frac{e_{c}}{r\beta})\alpha b_{\rm m}$, the co-evolutionary system can eventually reach the desired outcome, neither in the case of a moderately
growing resource pool nor the case of a rapidly growing resource pool. However, when the institutional inspection probability during a time unit becomes smaller such that $\hat{p}<(\frac{1}{\beta}-\frac{e_{c}}{r\beta})\alpha b_{\rm m}$, if no intervention is taken, the co-evolutionary system will evolve to a state with less resources and fewer cooperators and cannot eventually reach the desired outcome.

In order to make the feedback-evolving game model with the uncertain parameter $\hat{p}$ reach the desired outcome, in this paper we consider the nonlinear MRAC approach. This approach uses a stable reference model which can reach the desired outcome, and the objective is to design adaptive rules for specific control parameters such
that the system dynamics behave like the reference model. In particular, we choose the reference model which can produce the same nonlinear feedback-evolving game dynamics as the actual system given by,
\begin{equation}
\left\{
\begin{array}{ll}
\dot{x}_{\rm m}(t)=x_{\rm m}(t)[1-x_{\rm m}(t)][p^{*}\beta-\frac{\alpha b_{\rm m}y_{\rm m}(t)}{R_{\rm m}}]\\
\dot{y}_{\rm m}(t)=ry_{\rm m}(t)[1-\frac{y_{\rm m}(t)}{R_{\rm m}}]-N\frac{y_{\rm m}(t)}{R_{\rm m}}b_{\rm m}[1+(1-x_{\rm m})\alpha],\,\,\label{1.3}
\end{array}
\right.
\end{equation}
where $x_{\rm m}(t)$ denotes the expected fraction of cooperators in the population at time $t$ and $y_{\rm m}(t)$ denotes the expected amount of the common-pool resource at time $t$. In order to make the reference model reach the desired outcome as $t \rightarrow \infty$, here we consider the case where the resource pool does not grow slowly which satisfies $r>e_{c}$ and the institutional monitored intensity is not poorly which satisfies $p^{*}>(\frac{1}{\beta}-\frac{e_{c}}{r\beta})\alpha b_{\rm m}$.

Next, our goal is to design an adaptive law for adjusting the institutional inspection probability during a time unit $\hat{p}(t)$ in such a way that the state of the feedback-evolving game model with the uncertain parameter $\hat{p}$ given by Eq.~(\ref{1.2}) can track the state of the reference model given by Eq.~(\ref{1.3}). Finally, the feedback-evolving game model with the uncertain parameter $\hat{p}$ given by Eq.~(\ref{1.2}) can reach the desired outcome as $t \rightarrow \infty$.

For the sake of simplicity, we consider variables $s(t)=x(t)-1$, $v(t)=y(t)-(R_{\rm m}-\frac{Nb_{\rm m}}{r})$, $s_{\rm m}(t)=x_{\rm m}(t)-1$, and $v_{\rm m}(t)=y_{\rm m}(t)-(R_{\rm m}-\frac{Nb_{\rm m}}{r})$. Here $s(t)$ denotes the difference between the fractions of cooperators of the actual system and in the desired outcome,
$v(t)$ denotes the difference between the amounts of common-pool resources of the actual system and in the desired outcome,
$s_{\rm m}(t)$ denotes the difference between the fractions of cooperators of the reference model and in the desired outcome,
$v_{\rm m}(t)$ denotes the difference between the amounts of common-pool resources of the reference model and in the desired outcome.
Then the dynamical equations of the differences between the state of the actual system and the state of desired outcome can be written as
\begin{equation}
\left\{
\begin{array}{lll}
\dot{s}(t)&=(\alpha b_{\rm m}-\frac{\alpha Nb_{\rm m}^{2}}{R_{\rm m}r}-\hat{p}(t)\beta)s(t)+(\alpha b_{\rm m}-\frac{\alpha Nb_{\rm m}^{2}}{R_{\rm m}r}-\hat{p}(t)\beta)s^{2}(t)\\
            &+\frac{\alpha b_{\rm m}}{R_{\rm m}}s(t)v(t)+\frac{\alpha b_{\rm m}}{R_{\rm m}}s^{2}(t)v(t)\\
\dot{v}(t)&=(\alpha b_{\rm m}N-\frac{\alpha N^{2}b_{\rm m}^{2}}{R_{\rm m}r})s(t)+(\frac{Nb_{\rm m}}{R_{\rm m}}-r)v(t)\\
&-\frac{r}{R_{\rm m}}v^{2}(t)+\frac{\alpha Nb_{\rm m}}{R_{\rm m}}s(t)v(t)\,\,.\label{2}
\end{array}
\right.
\end{equation}
and the equations of the differences between the state of the reference model and the state of desired outcome can be written as
\begin{equation}
\left\{
\begin{array}{ll}
\dot{s}_{\rm m}(t)&=(\alpha b_{\rm m}-\frac{\alpha Nb_{\rm m}^{2}}{R_{\rm m}r}-p^{*}\beta)s_{\rm m}(t)+(\alpha b_{\rm m}-\frac{\alpha Nb_{\rm m}^{2}}{R_{\rm m}r}-p^{*}\beta)s^{2}_{\rm m}(t)\\
               &+\frac{\alpha b_{\rm m}}{R_{\rm m}}s_{\rm m}(t)v_{\rm m}(t)+\frac{\alpha b_{\rm m}}{R_{\rm m}}s^{2}_{\rm m}(t)v_{\rm m}(t)\\
\dot{v}_{\rm m}(t)&=(\alpha b_{\rm m}N-\frac{\alpha N^{2}b_{\rm m}^{2}}{R_{\rm m}r})s_{\rm m}(t)+(\frac{Nb_{\rm m}}{R_{\rm m}}-r)v_{\rm m}(t)\\
&-\frac{r}{R_{\rm m}}v^{2}_{\rm m}(t)+\frac{\alpha Nb_{\rm m}}{R_{\rm m}}s_{\rm m}(t)v_{\rm m}(t)
\,\,.\label{3}
\end{array}
\right.
\end{equation}

Since the state variables $x(t)$ and $x_{m}(t)$ both have the same desired value 1, and the state variables $y(t)$ and $y_{m}(t)$ both have the same desired value $(R_{\rm m}-\frac{Nb_{\rm m}}{r})$, then the object of adjusting $\hat{p}(t)$ for $x(t)\rightarrow x_{\rm m}(t)$ and $y(t)\rightarrow y_{\rm m}(t)$
as $t\rightarrow \infty$ is now changed to adjust $\hat{p}(t)$ for $s(t)\rightarrow s_{\rm m}(t)$ and $v(t)\rightarrow v_{\rm m}(t)$
as $t\rightarrow \infty$.

\begin{remark} \label{re1}
Next, we will design the update law for adjusting the institutional inspection probability during a time unit based on the Lyapunov
stability theory, such that $s(t)\rightarrow s_{\rm m}(t)$ and $v(t)\rightarrow v_{\rm m}(t)$
as $t\rightarrow \infty$.
\end{remark}

\section{Theoretical Results}

Define the error $e_{1}(t)=s(t)-s_{\rm m}(t)$ and $e_{2}(t)=v(t)-v_{\rm m}(t)$, then we have the error equation as

\begin{equation}
\left\{
\begin{array}{lll}
\dot{e}_{1}(t)=\left\{(\alpha b_{\rm m}-\frac{\alpha Nb_{\rm m}^{2}}{R_{\rm m}r}-p^{*}\beta)(1+2s_{\rm m}(t))+\frac{\alpha b_{\rm m}}{R_{\rm m}}(v_{\rm m}(t)+2s_{\rm m}(t)v_{\rm m}(t))\right\}e_{1}(t)\\
~~~~~~~~~+\frac{\alpha b_{\rm m}}{R_{\rm m}}(s_{\rm m}(t)+s_{\rm m}^{2}(t))e_{2}(t)+(\alpha b_{\rm m}-\frac{\alpha Nb_{\rm m}^{2}}{R_{\rm m}r}-p^{*}\beta)e_{1}^{2}(t)\\
~~~~~~~~~+\frac{\alpha b_{\rm m}}{R_{\rm m}}\left\{(1+2s_{\rm m}(t))e_{1}(t)e_{2}(t)+v_{\rm m}(t)e_{1}^{2}(t)+e_{1}^{2}(t)e_{2}(t)\right\}\\
~~~~~~~~~-\hat{p}(t)\beta s(t)+p^{*}\beta s(t)-\hat{p}(t)\beta s^{2}(t)+p^{*}\beta s^{2}(t)\\
\dot{e}_{2}(t)=\left[\alpha b_{\rm m}N-\frac{\alpha N^{2}b_{\rm m}^{2}}{R_{\rm m}r}+\frac{\alpha Nb_{\rm m}}{R_{\rm m}}v_{\rm m}(t)\right]e_{1}(t)\\
~~~~~~~~~+\left[\frac{Nb_{\rm m}}{R_{\rm m}}-r-\frac{2r}{R_{\rm m}}v_{\rm m}(t)+\frac{\alpha Nb_{\rm m}}{R_{\rm m}}s_{\rm m}(t)\right]e_{2}(t)\\
~~~~~~~~~-\frac{r}{R_{\rm m}}e_{2}^{2}(t)+\frac{\alpha Nb_{\rm m}}{R_{\rm m}}e_{1}(t)e_{2}(t)\,\,.\label{4}
\end{array}
\right.
\end{equation}

Denote $e(t)=[e_{1}(t),e_{2}(t)]^{T}$,
$A_{m}=\left[
\begin{array}{lll}
\alpha b_{\rm m}-\frac{\alpha Nb_{\rm m}^{2}}{R_{\rm m}r}-p^{*}\beta&0\\
\alpha b_{\rm m}N-\frac{\alpha N^{2}b_{\rm m}^{2}}{R_{\rm m}r}&\frac{Nb_{\rm m}}{R_{\rm m}}-r\\
\end{array}\right]$,
$B_{m}(t)=\left[
\begin{array}{lll}
2(\alpha b_{\rm m}-\frac{\alpha Nb_{\rm m}^{2}}{R_{\rm m}r}-p^{*}\beta)s_{\rm m}(t)+\frac{\alpha b_{\rm m}}{R_{\rm m}}v_{\rm m}(t)&\frac{\alpha b_{\rm m}}{R_{\rm m}}s_{m}(t)\\
\frac{\alpha b_{\rm m}N}{R_{\rm m}}v_{\rm m}(t)&-\frac{2r}{R_{\rm m}}v_{\rm m}(t)+\frac{\alpha b_mN}{R_{\rm m}}s_{\rm m}(t)\\
\end{array}\right]$,
$C_{m}(t)=\left[
\begin{array}{lll}
2\frac{\alpha b_{\rm m}}{R_{\rm m}}s_{\rm m}(t)v_{\rm m}(t)&\frac{\alpha b_{\rm m}}{R_{\rm m}}s_{\rm m}^{2}(t)\\
0&0\\
\end{array}\right]$,
and $B_{p}=\left[
\begin{array}{ll}
1\\
0\\
\end{array}\right]$,
then the error equation can be written as the vector form
\begin{equation}
\dot{e}(t)=[A_{\rm m}+B_{\rm m}(t)+C_{\rm m}(t)]e(t)-B_{\rm p}\tilde{p}(t)\beta[s(t)+s^{2}(t)]+f(t),
\label{5}
\end{equation}
where $\tilde{p}(t)=\hat{p}(t)-p^{*}$ is the control parameter error representing the difference between the institutional inspection probability during a time unit and its preset value,
$f(t)=\left[
\begin{array}{ll}
f_{1}(t)\\
f_{2}(t)\\
\end{array}\right]$,
$f_{1}(t)=(\alpha b_{\rm m}-\frac{\alpha Nb_{\rm m}^{2}}{R_{\rm m}r}-p^{*}\beta)e_{1}^{2}(t)+\frac{\alpha b_{\rm m}}{R_{\rm m}}\{[1+2s_{\rm m}(t)]e_{1}(t)e_{2}(t)+v_{\rm m}(t)e_{1}^{2}(t)+e_{1}^{2}(t)e_{2}(t)\}$, and
$f_{2}(t)=-\frac{r}{R_{\rm m}}e_{2}^{2}(t)+\frac{\alpha Nb_{\rm m}}{R_{\rm m}}e_{1}(t)e_{2}(t)$.

Now our goal is to adjust $\hat{p}(t)$ so that the error can converge to zero, which means that $s(t)\rightarrow s_{\rm m}(t)$ and $v(t)\rightarrow v_{\rm m}(t)$
as $t\rightarrow \infty$.


Next, we will give one of the main results of this paper in Lemma \ref{le1}.

\begin{lemma}\label{le1}
Consider that the dynamical equation of the difference between the state of the actual system and the state of desired outcome with the uncertain parameter $\hat{p}$ is described by Eq.~(\ref{2}), and the equation of the difference between the state of the reference model and the state of desired outcome is given by Eq.~(\ref{3}).

(1) If the update law of institutional monitoring probability during a time unit is designed by
\begin{equation}
\begin{array}{lll}
\dot{\hat{p}}(t)=ae^{T}(t)QB_{p}\beta s(t)[1+s(t)],
\label{5.2}
\end{array}
\end{equation}

(2) the error and control parameter error satisfy the initial condition
\begin{equation}
\begin{array}{lll}
e^{T}(t_{0})Qe(t_{0})+\frac{\tilde{p}^{2}}{a}(t_{0})<\lambda_{min}(Q)m^{2}, ~~~\text{for all~~} t>0,
\label{5.22}
\end{array}
\end{equation}
where $a>0$ and $Q=\left[
\begin{array}{lll}
q_{11}&q_{12}\\
q_{12}&q_{22}\\
\end{array}\right]$ is a positive-definite matrix satisfying the Lyapunov equation
\begin{equation}
   A_{m}^{T}Q+QA_{m}=-I_{2},
  \label{7}
\end{equation}

then the tracking error $e(t)$ in  Eq.~(\ref{5}) is asymptotically stable,
$s(t)\rightarrow s_{m}(t)$, and $v(t)\rightarrow v_{m}(t)$ as $t\rightarrow \infty$.
\end{lemma}

Proof. See Appendix for the details.

\begin{remark}\label{re2} Lemma 3.1 shows that if the updated rule of institutional monitoring probability during a time unit is designed by Eq.~(\ref{5.2}) and the error and parameter error satisfy the initial condition given by Eq.~(\ref{5.22}), then the error will tend to 0. Accordingly, $s(t)\rightarrow s_{\rm m}(t)$ and $v(t)\rightarrow v_{\rm m}(t)$ as $t\rightarrow \infty$.
\end{remark}

\begin{remark}\label{re3} The updated law of institutional monitoring probability during a time unit designed by Eq.~(\ref{5.2}) measures how much the institutional monitoring probability during a time unit can deviate from its preset value.
\end{remark}

We define the tracking error between the state of the feedback-evolving game with the uncertain parameter $\hat{p}$ described by Eq.~(\ref{1.2}) and the state of the reference model given by Eq.~(\ref{1.3}) as $e_{x,y}=[e_{x},e_{y}]^{T}$, where $e_{x}(t)=x(t)-x_{\rm m}(t)$ and $e_{y}(t)=y(t)-y_{\rm m}(t)$.
Since $x(t)=s(t)+1$, $y(t)=v(t)+R_{\rm m}-\frac{Nb_{\rm m}}{r}$, $x_{\rm m}(t)=s_{\rm m}(t)+1$, and $y_{\rm m}(t)=v_{\rm m}(t)+R_{\rm m}-\frac{Nb_{\rm m}}{r}$,
then $e_{x,y}=e(t)$. Therefore, from Lemma \ref{le1}, we can equivalently derive Theorem \ref{th1} as follows.

\begin{theorem}\label{th1} Consider the feedback-evolving game with the uncertain parameter $\hat{p}$ described by Eq.~(\ref{1.2}), and the reference model given by Eq.~(\ref{1.3}).
If the update law of institutional monitoring probability during a time unit is designed as
\begin{equation}
\begin{array}{lll}
\dot{\hat{p}}(t)=ae_{x,y}^{T}(t)QB_{p}\beta x(t)[x(t)-1]
\label{5.8}
\end{array}
\end{equation}
and the tracking error and control parameter error satisfy the initial condition given by Eq.~(\ref{5.22}),
then the tracking error $e_{x,y}(t)$ is asymptotically stable, $x(t)\rightarrow x_{\rm m}(t)$, and $y(t)\rightarrow y_{\rm m}(t)$.
\end{theorem}

\begin{remark}\label{re4} Theorem \ref{th1} shows that if the update law of institutional monitoring probability during a time unit is designed by Eq. (\ref{5.8}) and the tracking error and parameter error satisfy the initial condition given by Eq. (\ref{5.22}), then the tracking  error will tend to 0. Accordingly, $x(t)\rightarrow x_{\rm m}(t)$ and $y(t)\rightarrow y_{\rm m}(t)$ as $t\rightarrow \infty$. Then the feedback-evolving game with the uncertain parameter $\hat{p}$ can track the reference model asymptomatically.
Since the reference model can reach the desired outcome $t\rightarrow \infty$, the feedback-evolving game with uncertain parameter $\hat{p}$ described by Eq.~(\ref{1.2}) can reach the desired outcome $t\rightarrow \infty$.
\end{remark}

\section{Numerical Examples}
In this section, we will give three examples representing three different cases of $\hat{p}$ to illustrate the effectiveness of our theoretical analysis.

{\large{\bf{Example 1:}}}
Consider the case of the institutional monitoring probability during a time unit $\hat{p}$ changing from the preset value $p^{*}$ to a smaller value in a moderately growing resource pool, such that $\hat{p}<(\frac{1}{\beta}-\frac{e_{c}}{r\beta})\alpha b_{\rm m}$ and $e_{c}<r<e_{d}$, which corresponds to the middle yellow area of  Fig. 1.
To verify the validity of our results, in the actual co-evolutionary system, we set $\hat{p}= 0.09$ for $0\leq t<1000$, $\hat{p}= 0.07$ for $1000\leq t<3000$, and adaptively adjust $\hat{p}(t)$ for $t\geq3000$. Other parameters are set as $r = 0.6$, $\alpha = 0.5$, $p^{*}=0.09$, $\beta= 0.5$, $N = 100$, $R_{\rm m} = 100$, and $b_{\rm m} = 0.5$.
Then, the reference model can be written as
\begin{equation}
\left\{
\begin{array}{ll}
\dot{x}_{\rm m}(t)=x_{\rm m}(t)[1-x_{m}(t)][0.045-\frac{0.25y_{\rm m}(t)}{100}]\\
\dot{y}_{\rm m}(t)=0.6y_{\rm m}(t)[1-\frac{y_{\rm m}(t)}{100}]-\frac{y_{\rm m}(t)}{2}[1+0.5(1-x_{\rm m})]\,\,.\label{6.1}
\end{array}
\right.
\end{equation}
Our numerical results are shown in Fig. 2. Fig. 2(a) shows the time evolution of cooperation level and the states of resource. We can see that in the time domain $0\leq t<1000$, this co-evolutionary system can reach the desired outcome where all participants choose the cooperation strategy and the renewable resource can be maintained at a sustainable level. In the time domain $1000\leq t<3000$, the institutional monitoring probability during a time unit changes from the preset value $p^{*}=0.09$ to $\hat{p}= 0.07$, we can see that there are less resources and fewer cooperators in the co-evolutionary system. In order to make the state of the actual co-evolutionary system track the state of the reference model, so that the actual co-evolutionary system reaches the desired outcome again, in the time domain $t\geq3000$ we design the update rule of $\hat{p}$ as
\begin{equation}
\dot{\hat{p}}(t)=0.00001[3087347.58(x-x_{\rm m})+705.645(y-y_{\rm m})]x(x-1),
\label{6.2}
\end{equation}
then we can see that this co-evolutionary system can reach the desired outcome again, which is consistent with the theoretical results predicted by Theorem \ref{th1}.
Fig. 2(b) shows the time evolution of the tracking error in the time domain $t\geq3000$. We can see that the tracking error asymptotically converges to zero.
Fig. 2(c) shows the time evolution of the institutional monitoring probability during a time unit in the time domain $t\geq3000$, confirming that the probability value of detecting a defector during a time unit is kept in the range $[0,1]$, which is in line with the definition of the institutional monitoring probability during a time unit.

\begin{figure}
\centering
\includegraphics[width=4.5in]{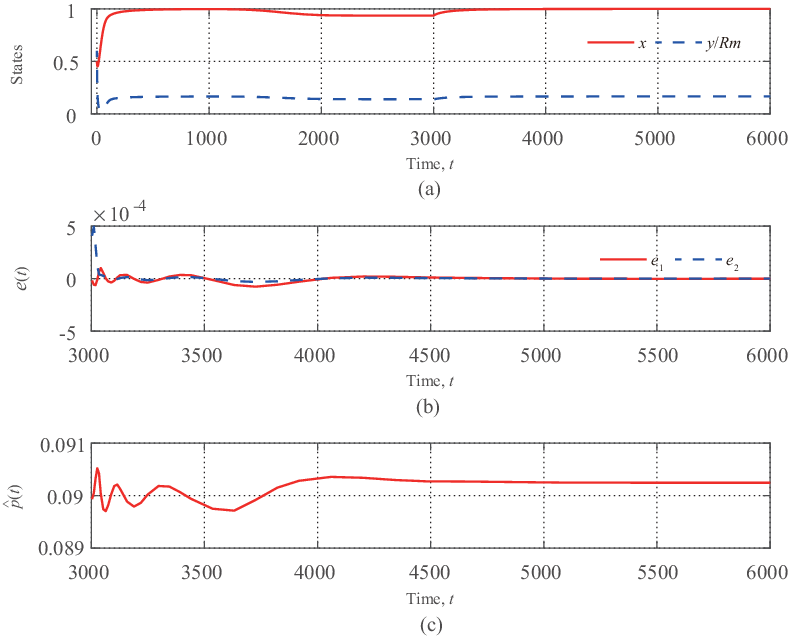}
\caption{Coevolutionary dynamics for $e_{c}<r < e_{d}$, $\hat{p}= 0.09$ for $0\leq t<1000$, $\hat{p}= 0.07$ for $1000\leq t<3000$, and adaptive $\hat{p}(t)$ for $t\geq3000$. Panel (a) shows the time evolution of cooperation level and the states of resource.
Panel (b) shows the time evolution of the tracking error in the time domain $t\geq3000$. Panel (c) shows the time evolution of the institutional monitoring probability during a time unit in the time domain $t\geq3000$. Other parameters are $r = 0.6$, $\alpha = 0.5$, $p^{*}= 0.09$, $\beta= 0.5$, $N = 100$, $R_{\rm m} = 100$, and $b_{\rm m} = 0.5$.
}
\label{fig2}
\end{figure}

{\large{\bf{Example 2:}}}
Consider the case of the institutional monitoring probability during a time unit $\hat{p}$ changing from the preset value $p^{*}$ to a smaller value in a
fast growing resource pool, such that $(\frac{1}{\beta}-\frac{e_{d}}{r\beta})\alpha b_{\rm m}<\hat{p}<(\frac{1}{\beta}-\frac{e_{c}}{r\beta})\alpha b_{\rm m}$ and $r>e_{d}$, which corresponds to the right yellow area of Fig. 1. To verify the validity of our results, in the actual co-evolutionary system, we set $\hat{p}= 0.2$ for $0\leq t<1000$, $\hat{p}= 0.17$ for $1000\leq t<4000$, and adaptive adjustly $\hat{p}(t)$ for $t\geq4000$. Other parameters are set as $r = 0.8$, $\alpha = 0.5$, $p^{*}=0.2$, $\beta= 0.5$, $N = 100$, $R_{\rm m} = 100$, and $b_{\rm m} = 0.5$.
Then, the expected reference model can be written as
\begin{equation}
\left\{
\begin{array}{ll}
\dot{x}_{\rm m}(t)=x_{\rm m}(t)[1-x_{\rm m}(t)][0.1-\frac{0.25y_{\rm m}(t)}{100}]\\
\dot{y}_{\rm m}(t)=0.8y_{\rm m}(t)[1-\frac{y_{\rm m}(t)}{100}]-\frac{y_{\rm m}(t)}{2}[1+0.5(1-x_{\rm m})]\,\,.\label{6.3}
\end{array}
\right.
\end{equation}
Fig. 3(a) shows the time evolution of cooperation level and the states of resource. We can see that in the time domain $0\leq t<1000$, this co-evolutionary system can reach the desired outcome. In the time domain $1000\leq t<4000$, the
institutional monitoring probability during a time unit changes from the preset value $p^{*}=0.2$ to $\hat{p}= 0.17$, we can see that there are less resources and fewer cooperators in the system. In order to make the state of the actual co-evolutionary system track the state of the reference model, so that the co-evolutionary system reaches the desired outcome again, in the time domain $t\geq4000$ we design the update law of $\hat{p}$ as
\begin{equation}
\dot{\hat{p}}(t)=0.0001[257307.89(x-x_{\rm m})+93.537(y-y_{\rm m})]x(x-1)
\label{6.4}
\end{equation}
then we can see that this actual co-evolutionary system can reach the desired outcome again, which is consistent with the theoretical results in Theorem \ref{th1}.
Fig. 3(b) shows the time evolution of the tracking error in the time domain $t\geq4000$. We can find that the tracking error asymptotically converges to zero.
Fig. 3(c) shows the time evolution of the institutional monitoring probability during a time unit in the time domain $t\geq4000$, confirming that the  probability value of detecting a defector during a time unit is kept in the range $[0,1]$, which is consistent with the definition of the institutional monitoring probability during a time unit.

\begin{figure}
\centering
\includegraphics[width=4.5in]{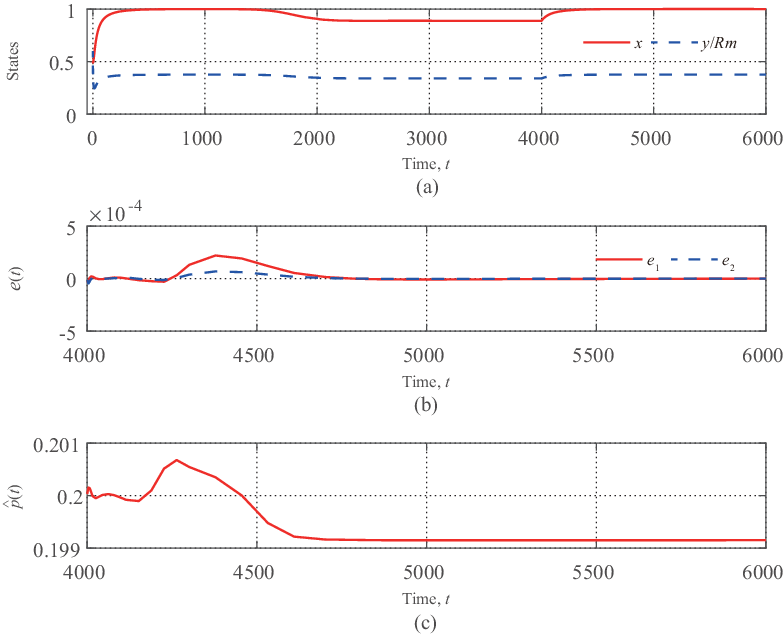}
\caption{Coevolutionary dynamics for $r > e_{d} $, $\hat{p}= 0.2$ for $t<1000$, $\hat{p}= 0.17$ for $1000\leq t<4000$, and adaptive $\hat{p}(t)$ for $t\geq4000$. Fig. 3(a) shows the time evolution of cooperation level and the states of resource.
Fig. 3(b) shows the time evolution of the tracking after $t>4000$.
Fig. 3(c) shows the time evolution of the institutional monitoring probability during a time unit. Other Parameters are $r = 0.8$, $\alpha = 0.5$, $p^{*}= 0.2$, $\beta= 0.5$, $N = 100$, $R_{\rm m} = 100$, and $b_{\rm m} = 0.5$.
}
\label{fig3}
\end{figure}

\begin{figure}
\centering
\includegraphics[width=4.5in]{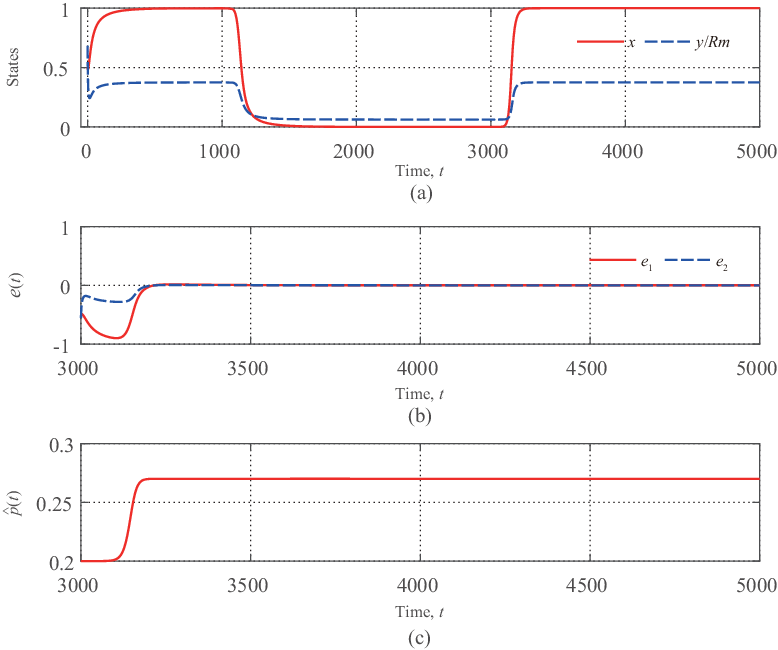}
\caption{Coevolutionary dynamics for $r > e_{d} $, $\hat{p}= 0.2$ for $t<1000$, $\hat{p}= 0.02$ for $1000\leq t<3000$, and adaptive $\hat{p}(t)$ for $t\geq3000$. Panel (a) shows the time evolution of cooperation level and the states of resource. Panel (b) shows the time evolution of the tracking error in the time domain $t>3000$. Panel (c) shows the time evolution of the institutional monitoring probability during a time unit in the time domain $t>3000$. Other parameters are $r = 0.8$, $\alpha = 0.5$, $p^{*}= 0.2$, $\beta= 0.5$, $N = 100$, $R_{\rm m} = 100$, and $b_{\rm m} = 0.5$.
}
\label{fig4}
\end{figure}

{\large{\bf{Example 3:}}}
Consider the case of the institutional monitoring probability during a time unit $\hat{p}$ changing from the preset value $p^{*}$ to a much smaller value in a
fast growing resource pool, such that $\hat{p}<(\frac{1}{\beta}-\frac{e_{d}}{r\beta})\alpha b_{\rm m}$ and $r>e_{d}$, which corresponds to the right red area of Fig.~1. To verify the validity of our results, in the actual co-evolutionary system we set $\hat{p}= 0.2$ for $0\leq t<1000$, $\hat{p}= 0.002$ for $1000\leq t<3000$, and adaptively adjust $\hat{p}(t)$ for $t\geq3000$. Other parameters are set as $r = 0.8$, $\alpha = 0.5$, $p^{*}=0.2$, $\beta= 0.5$, $N = 100$, $R_{\rm m} = 100$, and $b_{\rm m} = 0.5$.
Then, the expected reference model can be written as
\begin{equation}
\left\{
\begin{array}{ll}
\dot{x}_{\rm m}(t)=x_{\rm m}(t)[1-x_{\rm m}(t)][0.1-\frac{0.25y_{\rm m}(t)}{100}]\\
\dot{y}_{\rm m}(t)=0.8y_{\rm m}(t)[1-\frac{y_{\rm m}(t)}{100}]-\frac{y_{\rm m}(t)}{2}[1+0.5(1-x_{\rm m})]\,\,.\label{6.5}
\end{array}
\right.
\end{equation}
Fig. 4(a) shows the time evolution of cooperation level and the states of resource. We can see that in the time domain $0\leq t<1000$, this co-evolutionary system can reach the desired outcome. In the time domain $1000\leq t<3000$, the institutional monitoring probability during a time unit changes from the preset value $p^{*}=0.2$ to $\hat{p}= 0.002$, we can see that there are much fewer resources and no cooperators in the co-evolutionary system.
In order to drive the state of the actual co-evolutionary system to track the state of the reference model, so that the actual co-evolutionary system reaches the desired outcome again, in the time domain $t\geq3000$ we design the update law of $\hat{p}$ as
\begin{equation}
\dot{\hat{p}}(t)=0.0000001[257307.89(x-x_{\rm m})+93.537(y-y_{\rm m})]x(x-1)
\label{6.6}
\end{equation}
then we can see that this actual co-evolutionary system can reach the desired outcome again, which is consistent with the theoretical results in Theorem \ref{th1}.
Fig. 4(b) shows the time evolution of the tracking error in the time domain $t\geq3000$. We can observe that the tracking error asymptotically converges to zero.
Fig. 4(c) shows the time evolution of the institutional monitoring probability during a time unit in the time domain $t\geq3000$, confirming that the probability value of detecting a defector during a time unit is kept in the range $[0,1]$, which is in line with the definition of the institutional monitoring probability during a time unit.

\section{Conclusion}

In this work, we use the nonlinear MRAC approach to address the implementation uncertainty in institutional inspection intensity in a feedback-evolving game on the management of the common-pool renewable resource. We accordingly design a protocol that utilizes an update law to adjust the institutional inspection probability during a time unit based on the fraction of cooperators and the tracking error between the state of the actual system and the reference model. By means of Lyapunov stability theory, we derive a sufficient condition in which the update law can drive the actual system to reach the desired outcome in the feedback-evolving game system, regardless of the  growing rate of common-pool resource is moderate or rapid. Furthermore, we can measure how much the inspection probability can deviate from its preset value under this condition. In addition, we provide numerical examples for the moderately growing resource pool and fast growing resource pool respectively, in which the institutional inspection probability during a time unit is changed from the preset value to a smaller value, and all our numerical calculations in these examples verify our theoretical results. To our best knowledge, we use the nonlinear MRAC approach to address the uncertainty in a feedback-evolving game for the first time. This is a novel approach to address the implementation uncertainty in the coupled resource-strategy systems. We believe that our work can help to more effectively manage and restore the common-pool resources in human social systems.

\section*{Acknowledgments}
This research was supported by the National Natural Science
Foundation of China (Grants Nos. 61976048, 62036002, and 12171073).

\appendix

\section*{Appendix. The proof of Lemma 3.1}
\setcounter{section}{1}

Proof.
We choose the following quadratic Lyapunov function given by
\begin{equation}
	V(e(t),\tilde{p}(t))=e^{T}(t)Qe(t)+\frac{\tilde{p}^{2}}{a}(t).
	\label{8}
\end{equation}
Note that $V(0,0)=0$.
Since $Q$ is positive-definite, $V(e(t),\tilde{p}(t))>0$ for all $(e(t), \tilde{p}(t))\neq(0,0)$.

Differentiating Eq.~(\ref{8}) along the error equation described by Eq.~(\ref{5}) yields
\begin{equation}
	\begin{array}{lll}
		\dot{V}(e(t),\tilde{p}(t))
		&=\dot{e}^{T}(t)Qe(t)+e^{T}(t)Q\dot{e}(t)+\frac{2}{a}\tilde{p}(t)\dot{\hat{p}}(t)\\
		&=e^{T}(t)[A_{\rm m}+B_{\rm m}(t)+C_{\rm m}(t)]^{T}Qe(t)\\
		&~~~+e^{T}(t)Q[A_{\rm m}+B_{\rm m}(t)+C_{\rm m}(t)]e(t)\\
		&~~~-B_{\rm p}\beta\tilde{p}(t)[s(t)+s^{2}(t)]Qe(t)\\
		&~~~-e^{T}(t)QB_{\rm p}\beta\tilde{p}(t)[s(t)+s^{2}(t)]+\frac{2}{a}\tilde{p}\dot{\hat{p}}(t)\\
		&~~~+f^{T}(t)Qe(t)+e^{T}(t)Qf(t)\\
		&=e^{T}(t)(A_{\rm m}^{T}Q+QA_{m})e(t)\\
		&~~~+e^{T}(t)[(B_{\rm m}(t)+C_{\rm m}(t))^{T}Q+Q(B_{\rm m}(t)+C_{\rm m}(t))]e(t)\\
		&~~~-2e^{T}(t)QB_{\rm p}\beta\tilde{p}(s(t)+s^{2}(t))\\
		&~~~+\frac{2}{a}\tilde{p}\dot{\tilde{p}}(t)+2e^{T}(t)Qf(t)\\
		&=-e^{T}(t)e(t)+e^{T}(t)D(t)e(t)\\
		&~~~-2e^{T}(t)QB_{\rm p}\beta\tilde{p}(s(t)+s^{2}(t))\\
		&~~~+\frac{2}{a}\tilde{p}\dot{\hat{p}}(t)+2e^{T}(t)Qf(t),\\
	\end{array}
	\label{10}
\end{equation}
where $D(t)=(B_{\rm m}(t)+C_{\rm m}(t))^{T}Q+Q(B_{\rm m}(t)+C_{\rm m}(t))=\left[
\begin{array}{lll}
	d_{11}(t)&d_{12}(t)\\
	d_{12}(t)&d_{22}(t)\\
\end{array}\right]$,
with
$d_{11}(t)=4q_{11}(\alpha b_{\rm m}-\frac{\alpha Nb_{\rm m}^{2}}{R_{\rm m}r}-p^{*}\beta)s_{\rm m}(t)+(2q_{11}\frac{\alpha b_{\rm m}}{R_{\rm m}}+2q_{12}\frac{\alpha Nb_{\rm m}}{R_{\rm m}})v_{\rm m}(t)+4q_{11}\frac{\alpha b_{\rm m}}{R_{\rm m}}s_{\rm m}(t)v_{\rm m}(t)$,
$d_{12}(t)=[2q_{12}(\alpha b_{\rm m}-\frac{\alpha Nb_{\rm m}^{2}}{R_{\rm m}r}-p^{*}\beta)+q_{11}\frac{\alpha b_{\rm m}}{R_{\rm m}}+q_{12}\frac{\alpha Nb_{\rm m}}{R_{\rm m}}]s_{\rm m}(t)+(q_{12}\frac{\alpha b_{\rm m}}{R_{\rm m}}+q_{22}\frac{\alpha Nb_{\rm m}}{R_{\rm m}}-q_{12}\frac{2r}{R_{\rm m}})v_{\rm m}(t)+2q_{12}\frac{\alpha b_{\rm m}}{R_{\rm m}}s_{\rm m}(t)v_{\rm m}(t)+q_{11}\frac{\alpha b_{\rm m}}{R_{\rm m}}s_{\rm m}^{2}(t)$, and
$d_{22}(t)=(2q_{12}\frac{\alpha b_{\rm m}}{R_{\rm m}}+2q_{22}\frac{\alpha Nb_{\rm m}}{R_{\rm m}})s_{\rm m}(t)-q_{22}\frac{4r}{R_{\rm m}}v_{\rm m}(t)+2q_{12}\frac{\alpha b_{\rm m}}{R_{\rm m}}v_{\rm m}^{2}(t)$.

If the update law of institutional monitoring probability during a time unit is designed as Eq.~(\ref{5.2}),
then we can obtain that
\begin{equation}
	\begin{array}{lll}
		\dot{V}(e(t),\tilde{p}(t))&=-e^{T}(t)e(t)+e^{T}(t)D(t)e(t)+2e^{T}(t)Qf(t)\\
		&\leq-\| e(t)\|_{2}^{2}+\| e(t)\|_{2}\| D(t)\|_{F}\| e(t)\|_{2}+2\| e(t)\|_{2}\| Q\|_{F}\| f(t)\|_{2}\\
		&\leq-\| e(t)\|_{2}^{2}+\sqrt{d_{11}^{2}(t)+2d_{12}^{2}(t)+d_{22}^{2}(t)}\|e(t)\|_{2}^{2}\\
		&~~~+2\| e(t)\|_{2}\| Q\|_{F}\sqrt{f_{1}^{2}+f_{2}^{2}}.
	\end{array}
	\label{11}
\end{equation}

Note that Eq.~(\ref{3}) derived from the reference model Eq.~(\ref{1.3}) is asymptotically stable, so $s_{\rm m}(t)\rightarrow 0$ and $v_{\rm m}(t)\rightarrow 0$
as $t\rightarrow \infty$.
Therefore, $\forall \varepsilon>0$, there exists $T>0$, s.t. $\forall t>T$, $|s_{\rm m}(t)|<\varepsilon$ and $|v_{\rm m}(t)|<\varepsilon$.
We choose $\varepsilon<1$, then $\forall t>T$ we have $\varepsilon^{2}<\varepsilon$. Since $Q$ is positive definite, then $q_{11}>0$ and $q_{22}>0$.
Accordingly,
\begin{equation}
	\begin{array}{lll}
		|d_{11}(t)|&\leq4q_{11}|\alpha b_{\rm m}-\frac{\alpha Nb_{\rm m}^{2}}{R_{\rm m}r}-p^{*}\beta||s_{\rm m}(t)|\\
		&~~~~+(2q_{11}\frac{\alpha b_{\rm m}}{R_{\rm m}}+2|q_{12}|\frac{\alpha Nb_{\rm m}}{R_{\rm m}})|v_{\rm m}(t)|\\
		&~~~~+4q_{11}\frac{\alpha b_{\rm m}}{R_{\rm m}}|s_{\rm m}(t)||v_{\rm m}(t)|\\
		&<[4q_{11}|\alpha b_{\rm m}-\frac{\alpha Nb_{\rm m}^{2}}{R_{\rm m}r}-p^{*}\beta|+6q_{11}\frac{\alpha b_{\rm m}}{R_{\rm m}}+2|q_{12}|\frac{\alpha Nb_{\rm m}}{R_{\rm m}}]\varepsilon,
	\end{array}
	\label{11.2}
\end{equation}

\begin{equation}
	\begin{array}{lll}
		|d_{12}(t)|\leq[2|q_{12}(\alpha b_{\rm m}-\frac{\alpha Nb_{\rm m}^{2}}{R_{\rm m}r}-p^{*}\beta)|
		+q_{11}\frac{\alpha b_{\rm m}}{R_{\rm m}}+|q_{12}|\frac{\alpha Nb_{\rm m}}{R_{\rm m}}]|s_{\rm m}(t)|\\
		~~~~~~~~~~~~+(|q_{12}|\frac{\alpha b_{\rm m}}{R_{\rm m}}+q_{22}\frac{\alpha Nb_{\rm m}}{R_{\rm m}}+|q_{12}|\frac{2r}{R_{\rm m}})|v_{\rm m}(t)|\\
		~~~~~~~~~~~~+2|q_{12}|\frac{\alpha b_{\rm m}}{R_{\rm m}}|s_{\rm m}(t)||v_{\rm m}(t)|+q_{11}\frac{\alpha b_{\rm m}}{R_{\rm m}}|s_{\rm m}^{2}(t)|\\
		~~~~~~~~~~<[2|q_{12}(\alpha b_{\rm m}-\frac{\alpha Nb_{\rm m}^{2}}{R_{\rm m}r}-p^{*}\beta)|\\
		~~~~~~~~~~~~+(2q_{11}+N|q_{12}|+2|q_{12}|+Nq_{22})\frac{\alpha b_{\rm m}}{R_{\rm m}}\\
		~~~~~~~~~~~~~~+|q_{12}|\frac{2r}{R_{\rm m}}]\varepsilon,
	\end{array}
	\label{11.3}
\end{equation}
and
\begin{equation}
	\begin{array}{lll}
		|d_{22}(t)|&\leq(2|q_{12}|\frac{\alpha b_{\rm m}}{R_{\rm m}}+2q_{22}\frac{\alpha Nb_{\rm m}}{R_{\rm m}})|s_{\rm m}(t)|\\
		&~~~~+q_{22}\frac{4r}{R_{\rm m}}|v_{\rm m}(t)|
		+2|q_{12}|\frac{\alpha b_{\rm m}}{R_{\rm m}}|s_{\rm m}^{2}(t)|\\
		&<[(4|q_{12}|+2Nq_{22})\frac{\alpha b_{\rm m}}{R_{\rm m}}+q_{22}\frac{4r}{R_{\rm m}}]\varepsilon,
	\end{array}
	\label{11.4}
\end{equation}

\begin{equation}
	\begin{array}{lll}
		\|D(t)\|_{F}&=\sqrt{d_{11}^{2}(t)+2d_{12}^{2}(t)+d_{22}^{2}(t)}\\
		&\leq k\varepsilon,
	\end{array}
	\label{11.5}
\end{equation}
where
$k=\sqrt{k_{1}^2+2k_{2}^2+k_{3}^2}$
with
$k_{1}=4q_{11}|\alpha b_{\rm m}-\frac{\alpha Nb_{\rm m}^{2}}{R_{\rm m}r}-p^{*}\beta|+6q_{11}\frac{\alpha b_{\rm m}}{R_{\rm m}}+2|q_{12}|\frac{\alpha Nb_{\rm m}}{R_{\rm m}}$,
$k_{2}=2|q_{12}(\alpha b_{\rm m}-\frac{\alpha Nb_{\rm m}^{2}}{R_{\rm m}r}-p^{*}\beta)|+(2q_{11}+N|q_{12}|+2|q_{12}|+Nq_{22})\frac{\alpha b_{\rm m}}{R_{\rm m}}+|q_{12}|\frac{2r}{R_{\rm m}}$,
$k_{3}=(4|q_{12}|+2Nq_{22})\frac{\alpha b_{\rm m}}{R_{\rm m}}+q_{22}\frac{4r}{R_{\rm m}}$,
\begin{equation}
	\begin{array}{lll}
		f_{1}^{2}(t)&=(\alpha b_{\rm m}-\frac{\alpha Nb_{\rm m}^{2}}{R_{\rm m}r}-p^{*}\beta)^{2}e_{1}^{4}(t)\\
		&+\frac{\alpha^{2}b_{\rm m}^{2}}{R_{\rm m}^{2}}\{[1+4s_{\rm m}(t)+4s_{\rm m}^{2}(t)]e_{1}^{2}(t)e_{2}^{2}(t)+v_{\rm m}^2(t) e_{1}^{4}(t)+e_{1}^{4}(t)e_{2}^{2}(t)\}\\
		&+2(\alpha b_{\rm m}-\frac{\alpha Nb_{\rm m}^{2}}{R_{\rm m}r}-p^{*}\beta)\{\frac{\alpha b_{m}}{R_{\rm m}}
		[(1+2s_{\rm m}(t))e_{1}^{3}(t)e_{2}(t)+v_{\rm m}(t)e_{1}^{4}(t)]\\
		&+e_{1}^{4}(t)e_{2}(t)\}
		+2\frac{\alpha^{2}b_{\rm m}^{2}}{R_{\rm m}^{2}}[1+2s_{\rm m}(t)][v_{\rm m}(t)e_{1}^{3}(t)e_{2}(t)+e_{1}^{3}(t)e_{2}^{2}(t)]\\
		&+2\frac{\alpha^{2}b_{\rm m}^{2}}{R_{\rm m}^{2}}v_{\rm m}(t) e_{1}^{4}(t)e_{2}(t)\\
		&\leq(\alpha b_{\rm m}-\frac{\alpha Nb_{\rm m}^{2}}{R_{\rm m}r}-p^{*}\beta)^{2}e_{1}^{4}(t)\\
		&+\frac{\alpha^{2}b_{\rm m}^{2}}{R_{\rm m}^{2}}[(1+4|s_{\rm m}(t)|+4|s_{\rm m}^{2}(t)|)e_{1}^{2}(t)e_{2}^{2}(t)
		+|v_{\rm m}^2(t)|e_{1}^{4}(t)+e_{1}^{4}(t)e_{2}^{2}(t)]\\
		&+2|\alpha b_{\rm m}-\frac{\alpha Nb_{\rm m}^{2}}{R_{\rm m}r}-p^{*}\beta|\{\frac{\alpha b_{\rm m}}{R_{\rm m}}[(1+2|s_{m}(t)|)|e_{1}^{3}(t)e_{2}(t)|+|v_{\rm m}(t)|e_{1}^{4}(t)]\\
		&+e_{1}^{4}(t)|e_{2}(t)|\}+2\frac{\alpha^{2}b_{\rm m}^{2}}{R_{\rm m}^{2}}(1+2|s_{\rm m}(t)|)[|v_{\rm m}(t)||e_{1}(t)|^{3}|e_{2}(t)|+|e_{1}(t)|^{3}e_{2}^{2}(t)]\\
		&+2\frac{\alpha^{2}b_{\rm m}^{2}}{R_{\rm m}^{2}}|v_{\rm m}(t)|e_{1}^{4}(t)|e_{2}(t)|\\
		&\leq(\alpha b_{\rm m}-\frac{\alpha Nb_{\rm m}^{2}}{R_{\rm m}r}-p^{*}\beta)^{2}e_{1}^{4}(t)\\
		&+\frac{\alpha^{2}b_{\rm m}^{2}}{R_{\rm m}^{2}}[(1+4\varepsilon+4\varepsilon)e_{1}^{2}(t)e_{2}^{2}(t)+\varepsilon e_{1}^{4}(t)+e_{1}^{4}(t)e_{2}^{2}(t)]\\
		&+2|\alpha b_{\rm m}-\frac{\alpha Nb_{\rm m}^{2}}{R_{\rm m}r}-p^{*}\beta|\{{\frac{\alpha b_{\rm m}}{R_{\rm m}}[(1+2\varepsilon)|e_{1}^{3}(t)e_{2}(t)|+\varepsilon e_{1}^{4}(t)]+e_{1}^{4}(t)|e_{2}(t)|}\}\\
		&+2\frac{\alpha^{2}b_{\rm m}^{2}}{R_{\rm m}^{2}}(1+2\varepsilon)[\varepsilon|e_{1}^{3}(t)e_{2}(t)|+|e_{1}^{3}(t)|e_{2}^{2}(t)]\\
		&+2\frac{\alpha^{2}b_{\rm m}^{2}}{R_{\rm m}^{2}}\varepsilon e_{1}^{4}(t)|e_{2}(t)|,\\
	\end{array}
	\label{12}
\end{equation}
and
\begin{equation}
	\begin{array}{lll}
		f_{2}^{2}(t)&=[-\frac{r}{R_{\rm m}}e_{2}^{2}(t)+\frac{\alpha Nb_{\rm m}}{R_{\rm m}}e_{1}(t)e_{2}(t)]^{2}\\
		&=\frac{r^{2}}{R_{\rm m}^{2}}e_{2}^{4}(t)+\frac{\alpha^{2}N^{2}b_{\rm m}^{2}}{R_{\rm m}^{2}}e_{1}^{2}(t)e_{2}^{2}(t)-2\frac{r}{R_{\rm m}}\frac{\alpha Nb_{\rm m}}{R_{\rm m}}e_{1}(t)e_{2}^{3}(t).\\
	\end{array}
	\label{13}
\end{equation}
For the case of $e_{1}(t)<e_{2}(t)$,
\begin{equation}
	\begin{array}{llll}
		f_{1}^{2}&\leq(\alpha b_{\rm m}-\frac{\alpha Nb_{\rm m}^{2}}{R_{\rm m}r}-p^{*}\beta)^{2}e_{2}^{4}(t)\\
		&+\frac{\alpha^{2}b_{\rm m}^{2}}{R_{\rm m}^{2}}[(1+4\varepsilon+4\varepsilon)e_{2}^{2}(t)e_{2}^{2}(t)+\varepsilon e_{2}^{4}(t)+e_{2}^{4}(t)e_{2}^{2}(t)]\\
		&+2|\alpha b_{\rm m}-\frac{\alpha Nb_{\rm m}^{2}}{R_{\rm m}r}-p^{*}\beta|({\frac{\alpha b_{\rm m}}{R_{\rm m}}[(1+2\varepsilon)|e_{2}^{3}(t)e_{2}(t)|+\varepsilon e_{2}^{4}(t)]+e_{2}^{4}(t)|e_{2}(t)|})\\
		&+2\frac{\alpha^{2}b_{\rm m}^{2}}{R_{\rm m}^{2}}[3\varepsilon|e_{2}^{3}(t)||e_{2}(t)|+(1+2\varepsilon)|e_{2}^{3}(t)|e_{2}^{2}(t)]\\
		&+2\frac{\alpha^{2}b_{\rm m}^{2}}{R_{\rm m}^{2}}\varepsilon e_{2}^{4}(t)|e_{2}(t)|\\
		&\leq(\alpha b_{\rm m}-\frac{\alpha Nb_{\rm m}^{2}}{R_{\rm m}r}-p^{*}\beta)^{2}e_{2}^{4}(t)
		+\frac{\alpha^{2}b_{\rm m}^{2}}{R_{\rm m}^{2}}[10\varepsilon e_{2}^{4}(t)+e_{2}^{6}(t)]\\
		&+2|\alpha b_{\rm m}-\frac{\alpha Nb_{\rm m}^{2}}{R_{\rm m}r}-p^{*}\beta|[\frac{\alpha b_{\rm m}}{R_{\rm m}}(1+3\varepsilon)e_{2}^{4}(t)+|e_{2}^{5}(t)|]\\
		&+2\frac{\alpha^{2}b_{\rm m}^{2}}{R_{\rm m}^{2}}[3\varepsilon e_{2}^{4}(t)+(1+2\varepsilon)|e_{2}^{5}(t)|]
		+2\frac{\alpha^{2}b_{\rm m}^{2}}{R_{\rm m}^{2}}\varepsilon |e_{2}^{5}(t)|.\\
	\end{array}
	\label{14}
\end{equation}

Here we suppose $\|e(t)\|_{2}<b$, since $|e_{2}(t)|\leq\|e(t)\|_{2}$, then we can obtain that
\begin{equation}
	\begin{array}{llll}
		f_{1}^{2}&\leq(\alpha b_{\rm m}-\frac{\alpha Nb_{\rm m}^{2}}{R_{\rm m}r}-p^{*}\beta)^{2}\|e(t)\|_{2}^{4}
		+\frac{\alpha^{2}b_{\rm m}^{2}}{R_{\rm m}^{2}}[10\varepsilon \|e(t)\|_{2}^{4}+\|e(t)\|_{2}^{6}(t)]\\
		&~~~~+2|\alpha b_{\rm m}-\frac{\alpha Nb_{\rm m}^{2}}{R_{\rm m}r}-p^{*}\beta|[\frac{\alpha b_{\rm m}}{R_{\rm m}}(1+3\varepsilon)\|e(t)\|_{2}^{4}(t)+\|e(t)\|_{2}^{5}]\\
		&~~~~+2\frac{\alpha^{2}b_{\rm m}^{2}}{R_{\rm m}^{2}}[3\varepsilon \|e(t)\|_{2}^{4}(t)+(1+2\varepsilon)\|e(t)\|_{2}^{5}]
		+2\frac{\alpha^{2}b_{\rm m}^{2}}{R_{\rm m}^{2}}\varepsilon \|e(t)\|_{2}^{5}\\
		&\leq[(\alpha b_{\rm m}-\frac{\alpha Nb_{\rm m}^{2}}{R_{\rm m}r}-p^{*}\beta)^{2}+2|\alpha b_{\rm m}-\frac{\alpha Nb_{\rm m}^{2}}{R_{\rm m}r}-p^{*}\beta|\frac{\alpha b_{m}}{R_{\rm m}}]\|e(t)\|_{2}^{4}\\
		&~~~~+\varepsilon(16\frac{\alpha^{2}b_{\rm m}^{2}}{R_{\rm m}^{2}}+6|\alpha b_{\rm m}-\frac{\alpha Nb_{\rm m}^{2}}{R_{\rm m}r}-p^{*}\beta|\frac{\alpha b_{\rm m}}{R_{\rm m}})\|e(t)\|_{2}^{4}\\
		&~~~~+[2|\alpha b_{\rm m}-\frac{\alpha Nb_{\rm m}^{2}}{R_{\rm m}r}-p^{*}\beta|+2\frac{\alpha^{2}b_{\rm m}^{2}}{R_{\rm m}^{2}}(1+3\varepsilon)]b\|e(t)\|_{2}^{4}
		+\frac{\alpha^{2}b_{\rm m}^{2}}{R_{\rm m}^{2}}b^{2}\|e(t)\|_{2}^{4}\\
		&\leq l_{1}\|e(t)\|_{2}^{4},
	\end{array}
	\label{14.1}
\end{equation}
where $l_{1}=(\alpha b_{\rm m}-\frac{\alpha Nb_{\rm m}^{2}}{R_{\rm m}r}-p^{*}\beta)^{2}
+2|\alpha b_{\rm m}-\frac{\alpha Nb_{\rm m}^{2}}{R_{\rm m}r}-p^{*}\beta|\frac{\alpha b_{m}}{R_{\rm m}}+\varepsilon(16\frac{\alpha^{2}b_{\rm m}^{2}}{R_{\rm m}^{2}}+6|\alpha b_{\rm m}-\frac{\alpha Nb_{\rm m}^{2}}{R_{\rm m}r}-p^{*}\beta|\frac{\alpha b_{\rm m}}{R_{\rm m}})+b[2|\alpha b_{\rm m}-\frac{\alpha Nb_{\rm m}^{2}}{R_{\rm m}r}-p^{*}\beta|+(2+6\varepsilon)\frac{\alpha^{2}b_{\rm m}^{2}}{R_{\rm m}^{2}}]+
b^{2}\frac{\alpha^{2}b_{\rm m}^{2}}{R_{\rm m}^{2}}$,
and
\begin{equation}
	\begin{array}{lll}
		f_{2}^{2}&\leq\frac{r^{2}}{R_{\rm m}^{2}}e_{2}^{4}(t)+\frac{\alpha^{2}N^{2}b_{\rm m}^{2}}{R_{\rm m}^{2}}e_{2}^{4}(t)-2\frac{r}{R_{\rm m}}\frac{\alpha Nb_{\rm m}}{R_{\rm m}}e_{1}(t)e_{2}^{3}(t)\\
		&\leq l_{2}\| e(t)\|_{2}^{4},
	\end{array}
	\label{15}
\end{equation}
where $l_{2}=(\frac{r}{R_{\rm m}}-\frac{\alpha Nb_{\rm m}}{R_{\rm m}})^{2}$.

For the case of $e_{1}(t)\geq e_{2}(t)$,
\begin{equation}
	\begin{array}{lll}
		f_{1}^{2}&\leq(\alpha b_{\rm m}-\frac{\alpha Nb_{\rm m}^{2}}{R_{\rm m}r}-p^{*}\beta)^{2}e_{1}^{4}(t)\\
		&+\frac{\alpha^{2}b_{\rm m}^{2}}{R_{\rm m}^{2}}[(1+4\varepsilon+4\varepsilon)e_{1}^{2}(t)e_{1}^{2}(t)+\varepsilon e_{1}^{4}(t)+e_{1}^{4}(t)e_{1}^{2}(t)]\\
		&+2|\alpha b_{\rm m}-\frac{\alpha Nb_{\rm m}^{2}}{R_{\rm m}r}-p^{*}\beta|({\frac{\alpha b_{\rm m}}{R_{\rm m}}[(1+2\varepsilon)|e_{1}^{3}(t)||e_{1}(t)|+\varepsilon e_{2}^{4}(t)]+e_{1}^{4}(t)|e_{1}(t)}|)\\
		&+2\frac{\alpha^{2}b_{\rm m}^{2}}{R_{\rm m}^{2}}[3\varepsilon|e_{1}^{3}(t)||e_{1}(t)|+(1+2\varepsilon)|e_{1}^{3}(t)|e_{1}^{2}(t)]\\
		&+2\frac{\alpha^{2}b_{\rm m}^{2}}{R_{\rm m}^{2}}\varepsilon e_{1}^{4}(t)|e_{1}(t)|\\
		&\leq(\alpha b_{\rm m}-\frac{\alpha Nb_{m}^{2}}{R_{\rm m}r}-p^{*}\beta)^{2}e_{1}^{4}(t)
		+\frac{\alpha^{2}b_{\rm m}^{2}}{R_{\rm m}^{2}}[10\varepsilon e_{1}^{4}(t)+e_{1}^{6}(t)]\\
		&+2|\alpha b_{\rm m}-\frac{\alpha Nb_{\rm m}^{2}}{R_{\rm m}r}-p^{*}\beta|[\frac{\alpha b_{\rm m}}{R_{\rm m}}(1+3\varepsilon)e_{1}^{4}(t)+|e_{1}^{5}(t)|]\\
		&+2\frac{\alpha^{2}b_{\rm m}^{2}}{R_{\rm m}^{2}}[3\varepsilon e_{1}^{4}(t)+(1+2\varepsilon)|e_{1}^{5}(t)|]
		+2\frac{\alpha^{2}b_{\rm m}^{2}}{R_{\rm m}^{2}}\varepsilon |e_{1}^{5}(t)|.\\
	\end{array}
	\label{16}
\end{equation}

Since $|e_{1}(t)|\leq\|e(t)\|_{2}$, we can obtain that
\begin{equation}
	\begin{array}{llll}
		f_{1}^{2}\leq l_{1}\|e(t)\|_{2}^{4}
	\end{array}
	\label{14.2}
\end{equation}
and
\begin{equation}
	\begin{array}{lll}
		f_{2}^{2}
		\leq l_{2}\| e(t)\|_{2}^{4}.
	\end{array}
	\label{17}
\end{equation}
Therefore,
\begin{equation}
	\|f(t)\|_{2}\leq K\|e(t)\|_{2}^{2},
	\label{18}
\end{equation}
where $K=\sqrt{l_{1}+l_{2}}$.

Substituting Eqs.~(\ref{11.5}) and (\ref{18}) into Eq.~(\ref{11}), we can obtain
\begin{equation}
	\begin{array}{lll}
		\dot{V}(e(t),\tilde{p}(t))&\leq-\| e(t)\|_{2}^{2}+k\varepsilon\|e(t)\|_{2}^{2}+ K\|e(t)\|_{2}\| Q\|_{F}\|e(t)\|_{2}^{2}\\
		&\leq-\| e(t)\|_{2}^{2}[1-k\varepsilon-2\|e(t)\|_{2}\| Q\|_{F}K].
	\end{array}
	\label{19}
\end{equation}

Let $m=\min\{\frac{1-k\varepsilon}{2\| Q\|_{F}K},b\}$, therefore, when $\| e(t)\|_{2}<m$, we have
\begin{equation}
	\dot{V}(e(t),\hat{p}(t))\leq0.
	\label{20}
\end{equation}
Since $\dot{V}\leq0$, $e(t)$ and $\hat{p}(t)$ are bounded.

In addition, it is shown that $\dot{V}(e(t),\tilde{p}(t))$ is uniformly continuous by examining the boundedness of its
derivative, where
\begin{equation}
	\begin{array}{lll}
		\ddot{V}(e(t),\tilde{p}(t))&=-2\dot{e}^{T}(t)e(t)+2\dot{e}^{T}(t)De(t)+2e^{T}(t)D\dot{e}(t)\\
		&~~~~~~+2\dot{e}^{T}(t)Qf(t)+2e^{T}(t)Q\dot{f}(t).
	\end{array}
	\label{4.9}
\end{equation}
Since $e(t)$ and $\tilde{p}(t)$ are bounded by the virtue of $\dot{V}(e(t),\tilde{p}(t))\leq0$,
$s(t)$, $v(t)$, $f(t)$, $\dot{e}(t)$, and $\dot{f}(t)$ are bounded because $e(t)$, $s_{\rm m}(t)$, $v_{\rm m}(t)$,  and $\hat{p}(t)$ are bounded, then $\ddot{V}(e(t),\tilde{p}(t))$ is bounded. Therefore, $\dot{V}(e(t),\tilde{p}(t))$
is uniformly continuous.
Now, we follow the Barbalat's lemma \cite{astrom_95} that when $\| e(t)\|_{2}<m$, $\dot{V}(e(t),\tilde{p}(t))\rightarrow0$,
and hence $e(t)\rightarrow0$ as $t\rightarrow \infty$. Therefore, we can conclude that the error is asymptotically stable.

Due to
\begin{equation}
	\lambda_{min}(Q)\|e(t)\|_{2}^{2}\leq \lambda_{min}(Q)\|e(t)\|_{2}^{2}+\frac{\tilde{p}(t)}{a}\leq V(e(t),\tilde{p}(t))< V(e(t_{0}),\tilde{p}(t_{0})),
	\label{21}
\end{equation}
we take $V(e(t_{0}),\tilde{p}(t_{0}))=c=\lambda_{min}(Q)m^{2}$, we can obtain that $\| e(t)\|_{2}<m$, then $\{(e(t_{0}),\tilde{p}(t_{0}))\in R^{^{3}}:V(e(t_{0}),\tilde{p}(t_{0}))< c\}$ is an estimate of the region of attraction. That is to say, when the initial condition of the error and the control parameter error satisfy Eq.~(\ref{5.22}), if the updated law of institutional monitoring probability during a time unit is designed by Eq.~(\ref{5.2}), the error will tend to 0. Accordingly, $s(t)\rightarrow s_{\rm m}(t)$ and $v(t)\rightarrow v_{\rm m}(t)$ as $t\rightarrow \infty$.

This completes the proof.


\section*{Reference}

\begin{thebibliography}{45}

\bibitem{ostrom_90}
Ostrom E. Governing the commons: the evolution of institutions for collective action. Cambridge, UK: Cambridge University Press 1990.

\bibitem{kraak_ff11}
Kraak SBM.
Exploring the `public goods game' model to overcome the tragedy of the commons in fisheries management.
Fish Fish 2011; 12: 18-33.

\bibitem{adams_03}
Adams WM, Brockington D, Dyson J, Vira B.
Managing tragedies: understanding conflict over common pool resources. Science 2003; 302: 1915-1916.

\bibitem{hardin_68}
Hardin G. The tragedy of the commons. Science 1968; 162: 1243-8.

\bibitem{hofbauer_98}
Hofbauer J, Sigmund K. Evolutionary games and population dynamics. Cambridge: Cambridge University Press; 1998.

\bibitem{smith_82}
Smith JM, Smith JMM. Evolution and the theory of games. Cambridge: Cambridge University Press; 1982.

\bibitem{weibull_97}
Weibull JW. Evolutionary game theory. Cambridge: MIT Press; 1997.

\bibitem{nowak_04}
Nowak MA, Sasaki A, Taylor C, Fudenberg D. Emergence of cooperation and evolutionary stability in finite populations. Nature 2004; 428: 646-50.

\bibitem{sigmund_n10}
Sigmund K, De Silva H, Traulsen A, Hauert C.
Social learning promotes institutions for governing the commons.
Nature 2010; 466: 861-863.

\bibitem{vasconcelos_ncc13}
Vasconcelos VV, Santos FC, Pacheco JM.
A bottom-up institutional approach to cooperative governance of risky commons.
Nat Clim Change 2013; 3: 797-801.

\bibitem{Chen-JRSI-15}
Chen X, Sasaki T,  Br{\"a}nnstr{\"o}m {\AA}, Dieckmann U. First carrot, then stick: how the adaptive hybridization of incentives promotes cooperation. J R Soc Interface 2015; 12: 20140935.

\bibitem{Chen-PRSB-19}
Chen X, Br{\"a}nnstr{\"o}m {\AA}, Dieckmann U. Parent-prefered dispersal promotes cooperation in structured populations. Proc R Soc B 2019; 286: 20181949.

\bibitem{Duong-PRSA-21}
Duong MH, Han TA. Cost efficiency of institutional incentives for promoting cooperation in finite populations. Proc R Soc A 2021; 477: 20210568.

\bibitem{Jusup-PR-22}
Jusup M, Holme P, Kanazawa K, Takayasu M, Romi{\'c} I, Wang Z, Ge{\v{c}}ek S, Lipi{\'c} T, Podobnik B, Wang L, Luo W, Klanj{\v{s}}{\v{c}}ek T, Fan J, Boccaletti S, Perc M. Social physics. Phys. Rep. 2022; 948: 1-148.

\bibitem{Quan-Chaos-23}
Quan J, Zhang X, Chen W, Wang X. Cooperation dynamics in collective risk games with endogenous endowments. Chaos 2023; 33: 073107.

\bibitem{Chen-ND-23}
Chen W, Quan J, Wang X, Liu Y. Evolutionary dynamics from fluctuating environments with deterministic and stochastic noises. Non Dyna 2023; 111: 5499-5511.

\bibitem{tavoni_jtb12}
Tavoni A, Schl\"{u}ter M, Levin S.
The survival of the conformist: social pressure and renewable resource management.
J Theor Biol 2012; 299: 152-61.

\bibitem{Chen-SR-14}
Chen X, Perc M. Excessive abundance of common resources deters social responsibility. Sci Rep 2014; 4: 4161.

\bibitem{weitz_pnas16}
Weitz JS, Eksin C, Paarporn K, Brown SP, Ratcliff WC.
An oscillating tragedy of the commons in replicator dynamics with game-environment feedback.
Proc Natl Acad Sci USA 2016; 113: E7518-E7525.

\bibitem{sugiarto_2017}
Sugiarto HS, Lansing JS, Chung NN, Lai CH, Cheong SA, Chew LY.
Social cooperation and disharmonary in communities mediated through common pool resource exploitation.
Phys Rev Lett 2017; 118: 208301.

\bibitem{chen_xj_pcb18}
Chen X, Szolnoki A.
Punishment and inspection for governing the commons in a feedback-evolving game.
PLoS Comput Biol 2018; 14: e1006347.

\bibitem{szolnoki_18}
Szolnoki A, Chen X. Environmental feedback drives cooperation in spatial social dilemmas. EPL 2018; 120: 58001.

\bibitem{shao_yx_epl19}
Shao Y, Wang X, Fu F.
Evolutionary dynamics of group cooperation with asymmetrical environmental feedback.
EPL 2019; 126: 40005.

\bibitem{hauert_jtb19}
Hauert C, Saade C, McAvoy A.
Asymmetric evolutionary games with environmental feedback.
J Theor Biol 2019; 462: 347-360.

\bibitem{lin_yh_prl19}
Lin YH, Weitz JS.
 Spatial interactions and oscillatory tragedies of the commons.
Phys Rev Lett 2019; 122: 148102.

\bibitem{Quan-JSM-19}
Quan J, Zhang M, Zhou Y, Wang X, Yang J-B. Dynamic scale return coefficient with environmental feedback promotes cooperation in spatial public goods game. J Stat Mech 2019: 103405.

\bibitem{wang_x_rspa20}
Wang X, Zheng Z, Fu F.
Steering eco-evolutionary game dynamics with manifold control.
Proc R Soc A 2020; 476: 20190643.

\bibitem{tilman_nc20}
Tilman AR, Plotkin JB, Ak\c{c}ay E.
Evolutionary games with environmental feedbacks.
Nat Commun 2020; 11: 915.

\bibitem{Hu-CSF-20}
Hu L, He N, Weng Q, Chen X, Perc M. Rewarding endowments lead to a win-win in the evolution of public cooperation and the accumulation of common resources. Chaos Soli Frac 2020; 134: 109694.

\bibitem{yan_21}
Yan F, Chen X, Qiu Z, Szolnoki A.
Cooperator driven oscillation in a time-delayed feedback-evolving game.
New J Phys 2021; 23: 053017.

\bibitem{cao_lx_21}
Cao LX, Wu B. Eco-evolutionary dynamics with payoff-dependent environmental
feedback. Chaos Soli Frac 2021; 150: 111088.

\bibitem{Shu-PRSA-22}
Shu L, Fu F. Eco-evolutionary dynamics of bimatrix games. Proc Roy Soc A 2022; 478: 20220567.

\bibitem{Ma-AMC-23}
Ma X, Quan J, Wang X. Evolution of cooperation with nonlinear environment feedback in repeated public goods game. Appl Math Comp 2023; 452: 128056.

\bibitem{Liu-elife-23}
Liu L, Chen X, Szolnoki A. Coevolutionary dynamics via adaptive feedback in collective-risk social dilemma game. eLife 2023; 12: e82954.

\bibitem{rice_96}
Rice JC, Richards LJ.
A framework for reducing implementation uncertainty in fisheries management. N Am J Fish Manage 1996; 16: 488-494.

\bibitem{dichmont_06}
Dichmont CM, Deng AR, Punt AE, Venables W,
Haddon M. Management strategies for short-lived species: the case of Australia's Northern Prawn
Fishery 1: accounting for multiple species, spatial
structure and implementation uncertainty when evaluating risk. Fish Res 2006; 82: 204-220.

\bibitem{fulton_11}
Fulton EA, Smith ADM, Smith DC,  Van Putten IE.
 Human behaviour: the key source of uncertainty in fisheries management.
 Fish Fish 2011; 12: 2-17.

\bibitem{arlinghaus_17}
Arlinghaus R, Al$\acute{o}$s J, Beardmore B, Daedlow K, Dorow M, Fujitani M, et al.
Understanding and managing freshwater recreational fisheries as complex adaptive social-ecological systems.
Rev Fish Sci Aquac 2017; 25: 1-41.

\bibitem{lavretsky09}
Lavretsky E.
Combined/Composite model reference adaptive control.
IEEE Trans Autom Control 2009; 54: 2692-2697.


\bibitem{yucelen13} Yucelen T,  Haddad WM.
Low-frequency learning and fast adaptation in model reference adaptive control.
IEEE Trans Autom Control 2013; 58: 1080-1085.

\bibitem{kersting18}
Kersting S, Buss M.
Direct and indirect model reference adaptive control for multivariable piecewise affine systems.
IEEE Trans Autom Control 2017; 62: 5634-5649.

\bibitem{oliveira18}
Oliveira TR, Rodrigues VHP, Fridman L.
Generalized model reference adaptive control by means of global HOSM differentiators.
IEEE Trans Autom Control 2019; 64: 2053-2060.

\bibitem{nguyen18}
Nguyen NT.
Model reference adaptive control:a primer.
Springer; 2018.

\bibitem{song19}
Song  G, Tao G.
A partial-state feedback model reference adaptive control scheme.
IEEE Trans Autom Control 2020; 65: 44-57.

\bibitem{qu19}
Qu Z, Thomsen BT, Annaswamy AM.
Adaptive control for a class of multi-input multi-output plants with arbitrary relative degree.
IEEE Trans Autom Control 2020; 65: 3023-3038.

\bibitem{yan_22}
Yan F,  Hou X,  Tian T.
Fractional order multivariable adaptive control based on a nonlinear scalar update law.
Mathematics 2022; 10: 3385.

\bibitem{yan_23}
Yan F, Hou X, Tian T. Model reference adaptive control based on a novel scalar update law. Int J Robust Nonlinear Control. 2023; 33: 6250-6262.

\bibitem{yermekbayeva_18}
Yermekbayeva JJ, Omarov AN, Nedorezov LV.
The research of adaptive control system for predator-prey model based on Lyapunov method.
In:2018 International Conference on Control, Artificial Intelligence, Robotics \& Optimization.
2018.

\bibitem{yuan_19}
Yuan C, Zeng W, Dai SL.
Distributed model reference adaptive containment control of heterogeneous uncertain multi-agent systems.
ISA Trans 2019; 86: 73-8.

\bibitem{nam88}
Nam K, Araposthathis A. A model reference adaptive control scheme for pure-feedback nonlinear systems. IEEE Trans Autom Control 1988; 33: 803-11.

\bibitem{scarritt08}
Scarritt S. Nonlinear model reference adaptive control for satellite attitude tracking. In AIAA Guidance, Navigation and Control Conference and Exhibit 2008.

\bibitem{asiain21}
Asiain E, Garrido R. Anti-chaos control of a servo system using nonlinear model reference adaptive control. Chaos Soli Frac 2021; 143: 110581.

\bibitem{astrom_95}
Astrom KJ, Wittenmark B. Adaptive Control. Springer; 1995.

\end{thebibliography}
\end{document}